\documentclass{article}
\usepackage{graphicx} 
\usepackage{amsfonts}
\usepackage{multirow}
\usepackage{tcolorbox}
\usepackage{amsmath}
\usepackage{amssymb}
\usepackage{tikz}
\usepackage{stmaryrd}
\usepackage{mathtools}
\usepackage{standalone}
\usetikzlibrary{arrows.meta, positioning, calc,shapes.geometric} 
\usepackage{capt-of}

\usepackage{url} 
\usepackage[utf8]{inputenc} 
\usepackage[T1]{fontenc}    
\usepackage{arxiv}
\usepackage[noblocks]{authblk}

\usepackage{hyperref}
\usepackage[backend=biber, style=alphabetic,doi=false,url=false]{biblatex}
\newcommand{\NN}{\mathbb{N}}
\newcommand{\RR}{\mathbb{R}}

\newcommand{\T}{\mathsf{T}}
\newcommand{\diag}{\operatorname{diag}}
\newcommand{\Nh}{N_{\mathrm{heads}}}
\newcommand{\dhead}{d_{\mathrm{head}}}

\title{Understanding Transformers and Attention Mechanisms: An Introduction for Applied Mathematicians}
\renewcommand{\shorttitle}{Understanding Transformers and Attention Mechanisms}

\author[1]{Michel Fabrice Serret}
\affil[1]{Center for Scientific Computing, Theory and Data, Paul Scherrer Institute, Switzerland.}
\affil[ ]{\texttt{\href{mailto:michel.serret@psi.ch}{michel.serret@psi.ch}}}
\date{}

\begin{document}
\maketitle
\begin{abstract}
     This document provides a brief introduction to the attention mechanism used in modern language models based on the Transformer architecture. We first illustrate how text is encoded as vectors and how the attention mechanism processes these vectors to encode semantic information. We then describe Multi-Headed Attention, examine how the Transformer architecture is built and look at some of its variants. Finally, we provide a glimpse at modern methods to reduce the computational and memory cost of attention, namely KV caching, Grouped Query attention and Latent Attention. This material is aimed at the applied mathematics community and was written as introductory presentation in the context of the IPAM Research Collaboration Workshop entitled \textit{“Randomized Numerical Linear Algebra” (RNLA)}, for the project: \textbf{Randomization in Transformer models}.
 
\end{abstract}

\section{Introduction}

\begin{figure}[h!]
\centering
    \includegraphics[width=\textwidth]{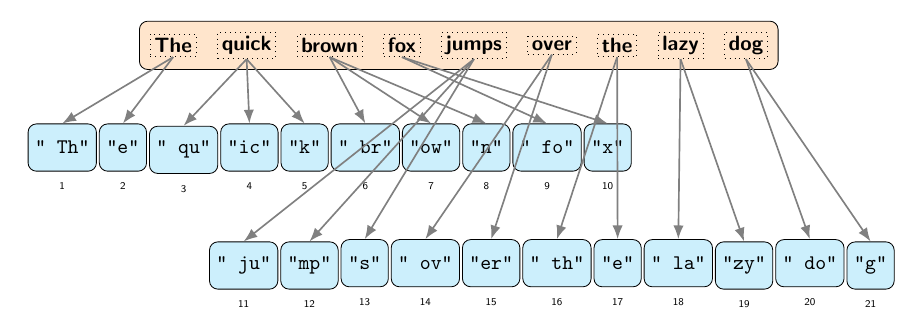}
\caption{Example of character tokenization of at most two characters.}
\label{fig:tokenization}
\end{figure}
Modern Natural Language Processing (NLP) methods based on attention mechanisms process textual information not in the form of strings of characters, but as sequences of vectors. 
To create these sequences, the text is divided into successive sub-strings called tokens, an example of which is given in Figure \ref{fig:tokenization}. The choice of the complete set of tokens used by the model, known as the vocabulary, is important.  Indeed, it needs to allow the encoding of enough of the semantic information contained in a text to achieve, given a model, the natural language processing tasks this model is designed to address. For example, the task consisting of predicting the next token in the text might require a finer tokenization than one classifying the complete sentence into two categories. However, it should also be as small as possible to avoid diluting the information contained in the sentence and to avoid unnecessarily increasing memory and computational costs.
The tokenization procedures vary from language to language and are often based on the structure of the language itself. Ideally, it would encode the smallest set of strings with inherent meaning, i.e. what would be defined as the set of all morphemes in linguistics. Second, once the text has been decomposed into a sequence of tokens, a vectorization procedure is required. In what follows, tokens can be taken without loss of generality to be words.

Now, given an ordered set $T=(T_i)_{1\le i\le N_{T}}$ of the $N_{T}$ possible distinct tokens, which we call the vocabulary, we can represent our text as a sequence of elements of $T$ through their index. This sequence of indices associated to a text is known as the tokenization. 

\begin{figure}[h!]
\centering
    \includegraphics[width=0.8\textwidth]{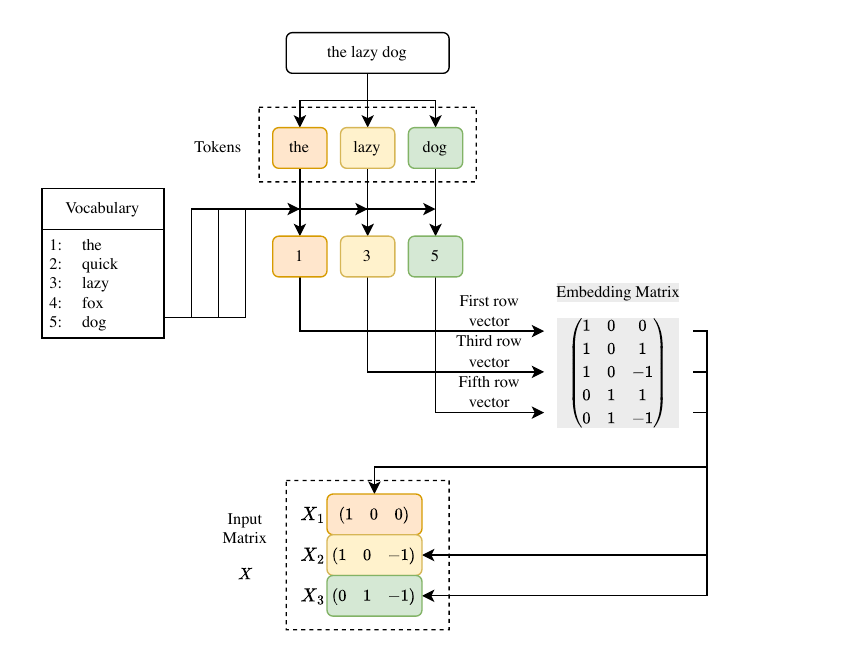}
\caption{Example of sentence embedding for a word-based tokenizer for the phrase ``the lazy dog''.}
\label{fig:tok_embd}
\end{figure}

Given a vocabulary, an embedding step allows to represent it in a lower dimensional space, through a trained linear application. This is also referred to as vectorization, since each token is thus represented by a vector.
Generally this embedding is given as a matrix  $E\in\RR^{N_{T}\times d}$ such that, given an embedding dimension $d\in\NN^*$, the $d$-dimensional vectorization of the $i$-th token $T_i$ in the vocabulary is given by the $i$-th row vector $E_{i}$ of $E$, as is illustrated in Figure \ref{fig:tok_embd}. These embedding matrices can either be pretrained or obtained through special methods and fine-tuned for the use cases of the language model or, as in most modern Large Language Models (LLMs), the embedding matrix can be trained from scratch as model weights.
To get an idea of the size of matrices at play, the vocabulary size of Llama 3 70b is $\sim 128$k with an embedding dimension $8192$ while for Gemma 3 27b it is $\sim 262$k of embedding dimension $5376$. Note that these are relatively small open source models.

Additionally, positional or feature embeddings (such as an embedding of the sentence to which the token belongs) can be applied to the sequence of token embedding vectors in order to enrich the semantic information contained therein. In some models, such as DeBERTa \cite{heDeBERTaDecodingenhancedBERT2021}or DeepSeek V2 \cite{deepseek-aiDeepSeekV2StrongEconomical2024}, the positional encoding and token embedding are treated separately, either through concatenation or by design.

\section{Attention mechanism}
\begin{figure}[h!]
\centering
    \includegraphics[width=\textwidth]{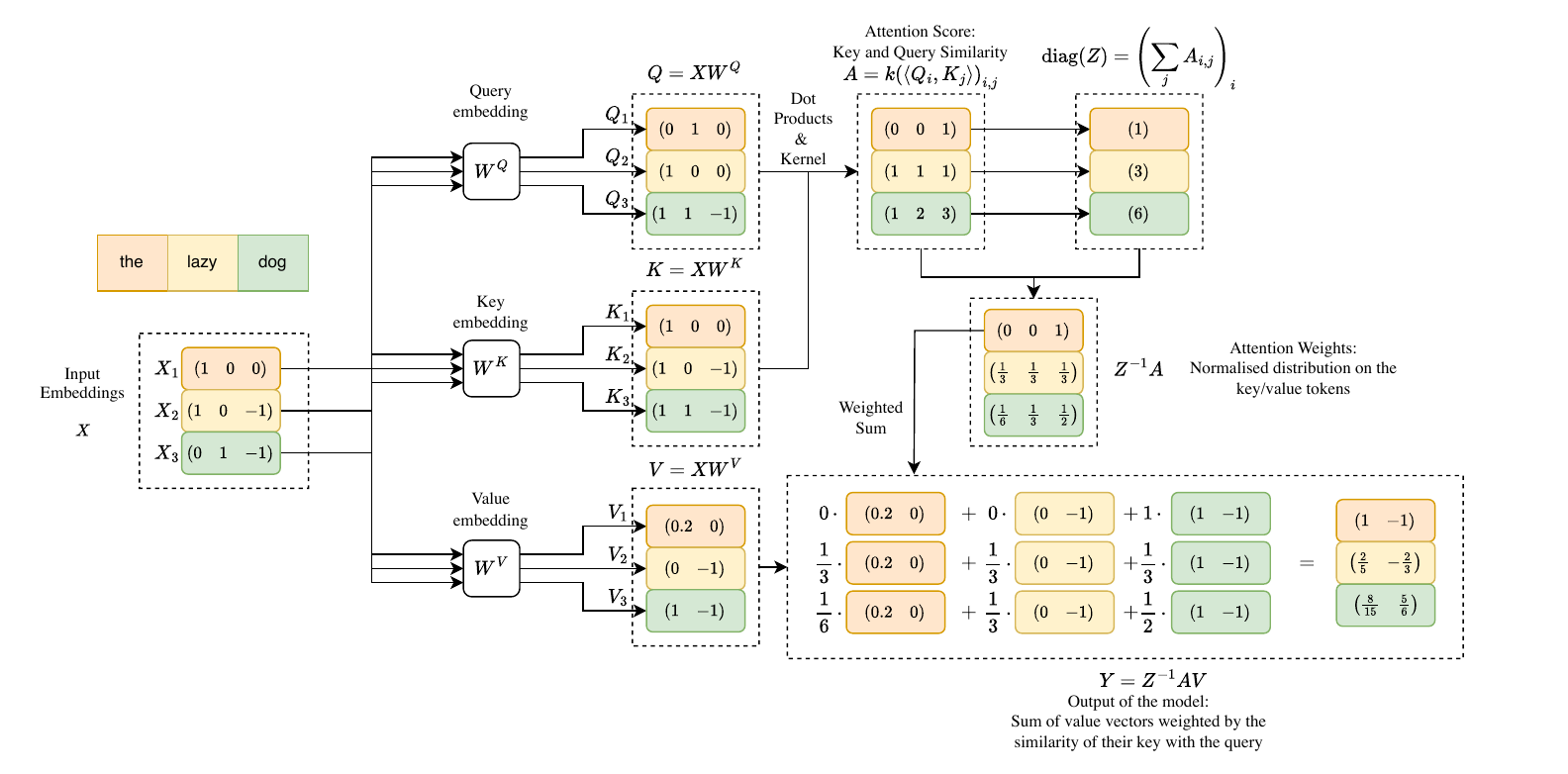}
\caption{Example of an attention layer for the embeddings obtained from the tokenization example in Figure \ref{fig:tok_embd}. Note that in this example we take the kernel to be linear for simplicity, i.e. for any $x\in\RR$, we set $f_\kappa(x)=x$.}
\label{fig:attn_meca}
\end{figure}

Attention mechanisms, the building blocks of the Transformer architecture \cite{vaswaniAttentionAllYou2017}, allow encoding of semantic information between tokens through a database-like structure\cite{zhangDiveDeepLearning2024}. 
Indeed, a database can be seen as a set of (key, value) tuples, $(k,v)$, such that when we submit a query $q$ to the database, the database returns the value $v$ associated to the tuple whose key is equal to the query, $q=k$. 
Analogously, let $X\in\RR^{N_{KV}\times d_{\text{in}}}$ be the set of $N_{KV}$ embedded tokens of dimension $d_{\text{in}}$ such as those obtained in Figure~\ref{fig:tok_embd}. As shown in Figure~\ref{fig:attn_meca}, given an input vector $x\in X$, associated to one of the tokens, the output is given as a linear combination of the ``value'' embeddings in the database, associated with each token $x'\in X$, weighted by a similarity metric relating the query embedding $q$ of $x$ to the ``key'' embedding $k'$ of $x'$.
As we shall see in the next section, in modern language models, attention mechanisms are used as modules, known in the literature as layers, and sequentially interweaved with other types of layers in order to construct the language model.
Let us now describe more formally the basic mechanism of a single ``layer'' of attention.
\subsection{General attention mechanism}
Consider the given matrix $X^Q\in\RR^{N_Q\times d_{\text{in}}^{Q}}$ of $N_Q\in\NN^*$ input query vectors, and the matrices $X^K\in\RR^{N_{KV}\times d_{\text{in}}^{K}}$  and $X^V\in\RR^{N_{KV}\times d_{\text{in}}^{V}}$  of  $N_{KV}\in\NN^*$ input key and value vectors respectively, with $d_{\text{in}}^{Q},d_{\text{in}}^{K},d_{\text{in}}^{V}\in\NN^*$ the respective dimensions of the vectors. In the following, we take  for simplicity, as the query, key and value inputs come from the same embeddings in most cases, $d_{\text{in}}^Q=d_{\text{in}}^K=d_{\text{in}}^V=d_{\text{in}}$. 
\begin{tcolorbox}
    \textbf{Example:} In the context of a conversation with chatGPT, the $N_Q$ query vectors, given by the rows of $Q$, would correspond to the $N_Q$ tokens in the last question asked to the model while the $N_{KV}$ key, encoded in the rows of $K$, and value vectors, encoded in the rows of $V$, would correspond to every token of the complete conversation. While in general the key and value vectors can be embedded from different sources for the same token, in the context of a chat, they would correspond to the same input vector, the embedding $X$ of the associated token as in the Figure~\ref{fig:tok_embd}.
\end{tcolorbox}
Furthermore, note that in general, the key and value inputs will be shared as they generally correspond to two embeddings of the same ``object'', while the query key can correspond to another object. As we shall see in the cross-attention module of the decoder of the original Transformer, the key and value vectors correspond to the tokens of a sentence in one language, while the query vectors correspond to the predicted tokens of the translation of the sentence in another language. For GPT models, such as chatGPT, the attention mechanism's inputs all come from the same source. The only difference in the query and key/value vectors arises from the fact that when a new question is asked to the model, only the new query vectors from the question need to be calculated if the key/value vectors from the discussion have been stored; this is known as streaming attention \cite{hanStreamingAttentionApproximation2025}.

The first step of the attention mechanism is to construct embeddings by applying a linear transformation to the query, key, and value input vectors. We denote the associated model weights as $W^{Q}\in\RR^{d_{\text{in}}^{Q}\times d_{QK}}$, $W^K\in\RR^{d_{\text{in}}^{K}\times d_{QK}}$ and $W^V\in\RR^{d_{\text{in}}^{V}\times d_{\text{out}}}$ where $d_{QK}$ and $d_{V}$ are the shared dimensions for the query and key embeddings and the value embedding dimension respectively, while these dimensions can be different in theory, in practice they are often equal to one another.
We denote the associated query, key and value matrices as
\[ 
    Q = X^QW^{Q}\in \RR^{N_{Q}\times d_{QK}},\quad K=X^KW^{K}\in \RR^{N_{KV}\times d_{QK}},\quad V=X^VW^{V}\in \RR^{N_{KV}\times d_{V}}.
\]
The similarity metric is usually given in the form of a kernel $\kappa:\RR^{d_{QK}}\times \RR^{d_{QK}}\to\RR$
which allows to obtain the \textbf{attention scores} for each tuple $(i,j)\in\llbracket1,N_{Q}\rrbracket\times \llbracket1,N_{KV}\rrbracket$ of the $N_Q$ query and $N_{KV}$ key/value indices,
$$A=\Big(\kappa(Q[i,:],K[j,:])\Big)_{\substack{1\le i\le N_{Q}\\1\le j\le N_{KV}}}\in\RR^{N_Q\times N_{KV}},$$
where we denote by $B[i,:]$ the $i$-th row vector of the matrix $B$.
Note that the kernel will often be of the form
\[
    \kappa(v,w)=f_{\kappa}\left(\langle v,w\rangle\right),
\] for a given function $f_{\kappa}:\RR\to\RR$ and we will use this notation throughout the document to simplify some equations. Furthermore, the element-wise application of $f_{\kappa}$ to a matrix $M$ will be denoted as $f_\kappa(M)$. 
As we wish to emulate a probability distribution on the set of key/value inputs for a given query, a normalization step is required; the normalization coefficients are obtained as the inverse of the diagonal matrix $Z\in\RR^{N_{Q}\times N_{Q}}$, whose non-zero elements are the sum of the rows of $A$, i.e. 
\[ 
    Z[i,i] =\sum\limits_{1\le j\le N_{KV}}\kappa\left(Q[i,:],K[j,:]\right),\quad\forall i\in\{1,\dots,N_{Q}\}.
\]
Normalized attention scores, i.e., the product $Z^{-1}A$, are known in the literature as \textbf{attention weights}.
The final output of the attention mechanism is then given by
$$Y=Z^{-1}AV\in\RR^{N_Q\times d_{\text{out}}}.$$
Note that the output dimension of the attention mechanism is given by the value dimension, i.e. $d_{\text{out}}=d_V$.
\begin{tcolorbox}
The most common kernel is the scaled exponential kernel used in the original Transformer \cite{vaswaniAttentionAllYou2017} given, for any $v,w\in\RR^{d_{QK}}$, by \[
\kappa(v,w)=\exp\left(\frac{\langle v,w\rangle}{\sqrt{d_{QK}}}\right).
\]

 It originates from the fact that if the entries of $v$ and $w$ come from an i.i.d. distribution with mean zero and variance $1$, then $\langle v,w\rangle$ has variance $d_{QK}$. Hence, the denominator rescales the scalar product so that its variance is 1. 
 Given the matrices $Q$ and $K$, the attention scores can be obtained as
 \[ 
     A=\exp\left(\frac{QK^\T}{\sqrt{d_{QK}}}\right),
 \]
 where the exponential is applied element-wise.
 The attention weights in this context are equal to the softmax function, $\sigma:\RR^N\to\RR^N$, applied to the query-key scalar product, the softmax function is defined as 
 \[ 
    \sigma :x\in\RR^{N}\mapsto \frac{1}{\sum\limits_{\ell=1}^{N}e^{x_\ell}}
    \begin{pmatrix}
    e^{x_{1}} & \dots & e^{x_N}
    \end{pmatrix}^\T,
 \] and hence, for any $i\in\{1,\dots ,N_{Q}\}$,
 \[(Z^{-1}A) [i,:] = \sigma\left(\frac{Q[i,:]K^\T}{\sqrt{d_{QK}}}\right).\]
\end{tcolorbox}

\subsection{Multi-headed attention}\label{sec:MHA}
In order to allow the attention mechanism to focus on more types of semantic information, multiple parallel attention heads can be used on the same input vectors.
The simplest way to do so is to consider stacked single attention heads with the constraint that they must all have the same output dimension $\dhead$, i.e. for any $h\in\llbracket1,\Nh\rrbracket$, $W^V_h\in\RR^{d_\text{in}^{V}\times \dhead}$.
Let $\Nh\in\NN^*$ be the number of attention heads, for any $h \in\llbracket1,\Nh\rrbracket$, we denote by $ W^{Q}_{h},W_{h}^{K}$ and $W_{h}^{V}$ the query, key and value embedding matrices  associated to the $h$-th head.
For each head $h\in\llbracket1,\Nh\rrbracket$, the attention score can then be obtained as
\[
A_h=\Big(\kappa\left({Q_{h}[i,:],K_{h}[j,:]}\right)\Big)_{\substack{1\le i\le N_{Q}\\1\le j\le N_{KV}}}=f_\kappa(Q_{h}K_{h}^{\T}),
\] where $Q_h=X^QW^Q_h\in\RR^{N_{Q}\times d_{QK}}$, $K_{h}={X^{K}}{W_{h}^{K}}\in\RR^{N_{KV}\times d_{QK}}$, and $k$ is applied element-wise. We note that in general $d_{QK}=\dhead$ for each head $h$.
Similarly, we set
\[ 
V_{h} =X^{V}W^{V}_{h},\quad Z_h=\diag\left(\sum\limits_{1\le j\le N_{KV}}\kappa\left({Q_{h}[i,:],K_{h}[j,:]}\right)\right)_{1\le i\le N_{Q}}.
\]
The output of each attention head can then be obtained as $$Y_{h}=Z_{h}^{-1}A_{h}V_{h}\in\RR^{N_{Q}\times \dhead}.$$
Given the outputs of each head, the output of the multi-headed attention mechanism is combined into a single output $Y\in\RR^{N_{Q}\times d_{\text{out}}}$ through a linear application $W^{O}\in\RR^{\Nh\dhead\times d_{\text{out}}}$, which is trained, as
\[ 
    Y=\operatorname{concat}_{\mathrm{col}}\left(Y_h\right)_{1\le h\le \Nh}W^O,
\]
where $\operatorname{concat}_{\mathrm{col}}\left(Y_h\right)_{1\le h\le \Nh}=\begin{pmatrix}
Y_1 &\cdots & Y_{\Nh}
\end{pmatrix}\in\RR^{N_{Q}\times \Nh \dhead}$.

\section{Transformer Architecture}
\subsection{Encoders, Decoders and the Transformer architecture}
In natural language processing, a common method used for machine translation of one language to another is the use of an encoder-decoder architecture to first encode an input sentence in the first language into an intermediary state and then use a decoder to transform this intermediary state into the equivalent sentence in the second language.
\begin{figure}[h!]\label{fig:enc_dec}
\centering
    \includegraphics[width=\textwidth]{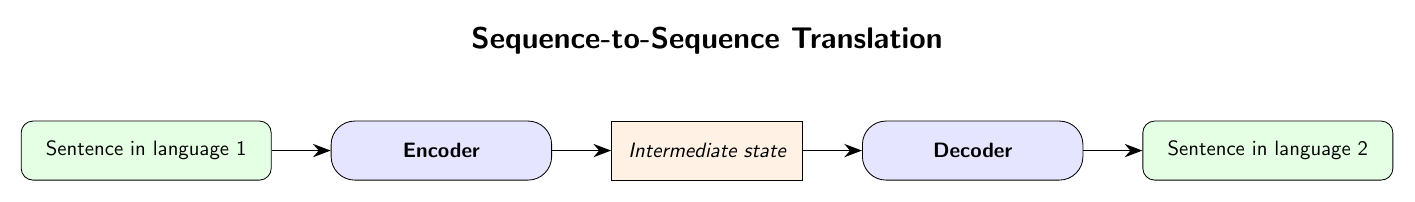}
\caption{Encoder-Decoder model}
\end{figure}

The original \textbf{Transformer} model \cite{vaswaniAttentionAllYou2017} was based on such a system. Let us now look at the individual components of the different layers which compose the Transformer architecture.
\subsubsection{Encoder layers}

\begin{figure}[h!]
    \centering
    \includegraphics[width=\linewidth]{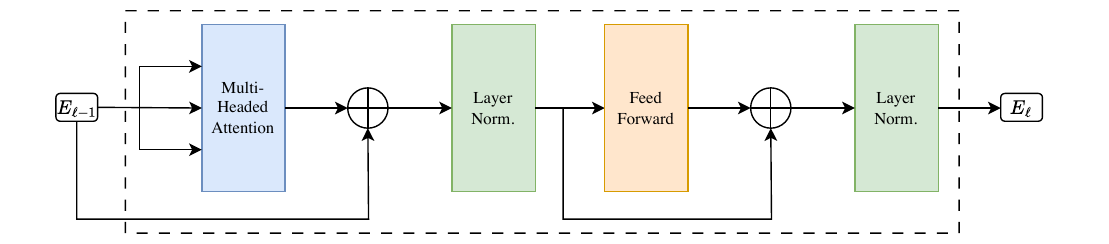}
    \caption{Encoder layer of the original Transformer architecture. Note that the $\bigoplus$ sign in the diagram corresponds to the skip connections.}
    \label{fig:t_enc}
\end{figure}
The encoder is composed of a series of $N_{L}$ layers of encoder layers . Each encoder layer consists of the following components, assembled as shown in Figure \ref{fig:t_enc}, 
\begin{enumerate}
    \item Multi-headed \textbf{self}-attention sublayer:\\ This sublayer consists in a multi-headed attention mechanism exactly as described in subsection \ref{sec:MHA}. The \textbf{self} part in the name refers to the fact that both its query and key/value input vectors are  from the same input, i.e. we have $X^Q=X^K=X^V$. In Figure~\ref{fig:t_enc} this corresponds to the fact that the three input arrows originate from the previous layer.
    \item Layer normalization sublayer: \\ 
		Given the input $x\in\RR^{d}$ to a layer, layer normalization\cite{baLayerNormalization2016} constructs a  vector $\tilde { x}=\frac{ x-\mu}{\sqrt{\sigma^2+\epsilon}}$ of mean $0$ and variance $1$, where $\mu=\frac{1}{d}\sum\limits_{i=1}^{d} x_{i}$ is the input vector's mean and $\sigma=\sqrt{\frac{1}{d}\sum_{i=1}^{d}(x_i-\mu)^2}$ is its standard deviation. It outputs  $ y=\gamma\odot \tilde{ x}+\beta$, where $\gamma,\beta\in\RR^{d}$ are known as the gain and shift and correspond to trainable parameters that rescale the variance and mean of the output to a fixed value and we denote by $\odot$ the Hadamard or term-wise product.
    \item  Feed forward sublayer: \\ 
		It is the basic building block of neural networks, it outputs $y=f(W x+b)$ where $f$ is a nonlinear function, $W\in\RR^{d\times d'}$ is the weight matrix of the sublayer and $b\in\RR^{d'}$ is the bias vector, where $d'$ is the dimension of the output of the layer. In general the nonlinear function $f$ is a variant of the Rectified Linear Unit (ReLU) family of functions; recently, parametric Gated Linear Unit (GLU) variants \cite{shazeerGLUVariantsImprove2020} have been favored in state-of-the-art models \cite{teamGemma2Improving2024,teamGemma3Technical2025,yangQwen3TechnicalReport2025,grattafioriLlama3Herd2024}. 
    \item Skip connections: \\
        Represented by the arrows that go under the components in Figure~\ref{fig:t_enc}, skip connections consist in adding to the output of a sublayer its input. They allow the preservation of the flow of information through the neural network and reduces the problem of vanishing gradients which arises in deep neural networks.
\end{enumerate}

\subsection{Decoder layers, Causal/Masked attention and Cross-attention}

\begin{figure}[h!]
    \centering
    \includegraphics[width=\linewidth]{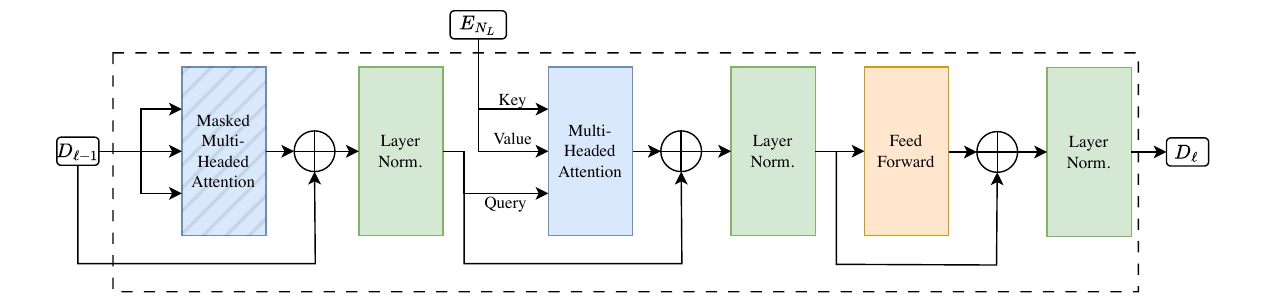}
    \caption{Decoder layer of the original Transformer architecture. Notice that the main differences with the encoder are the cross-attention sublayer where the key/value vectors are constructed not from the decoder state but from the final encoder state and the masked causal self-attention restricting the query vectors from seeing the key vectors of future tokens.}
    \label{fig:t_dec}   
\end{figure}
Decoder layers in the Transformer architecture have two differences with the encoder layers, as can be seen in Figure~\ref{fig:t_dec}.

\begin{tcolorbox}
        \textbf{Remark:}
        The architecture shown in Figure~\ref{fig:t_enc} is the original Transformer architecture. An alternative \cite{xiongLayerNormalizationTransformer2020}, reducing the problem of the vanishing gradients and allowing easier training but at the cost of performance, consists in applying the layer normalization before the attention sublayer and is shown in Figure~\ref{fig:PreLN}.
    Note how the skip connections are not normalized in this architecture and thus allow the gradients to skip through unscathed. Recently, state-of-the-art large language models such as Gemma 3\cites{teamGemma3Technical2025}, have made use of both Pre-Layer and Post-Layer Normalization\cite{kimPeriLNRevisitingLayer2025}. Furthermore, layer-normalization has often been replaced with the less computationally intensive RMS-Normalization procedure which consists in only rescaling the input with respect to its Root-mean-squared norm without recentering or adding a bias \cite{zhangRootMeanSquare2019}.
    
    {\centering
    \includegraphics[width=\linewidth]{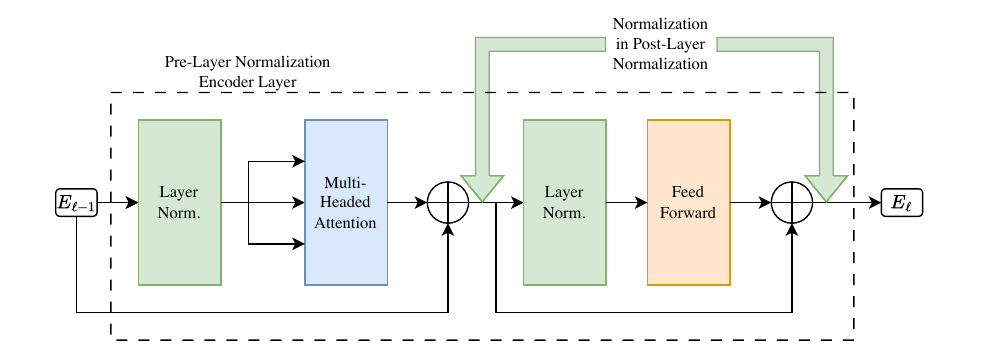}
    \captionof{figure}{Encoder layer with Pre-Layer Normalization}\label{fig:PreLN}}

\end{tcolorbox}
\paragraph{Masked and causal attention}
 The first difference is the use of masked attention. As the decoders construct the new sentence token by token, the self-attention mechanism that produces the token at step $t$ is not able to take into account the key-value vectors of tokens that will be created afterwards, hence, for consistency, a mask is applied to the attention to limit its vision of the tokens that come afterwards.
More formally, masked attention consists in restricting the subset of key-value tokens seen by each query token, i.e. defining, for each $1\le i\le N_Q$ , a subset of key-value tokens $S_i\subset\llbracket1,N_{KV} \rrbracket$ such that the attention matrix is given by
\[A=\Big(\mathbf 1_{j\in S_i}\kappa\left(Q[i,:],K[j,:]\right)\Big)_{\substack{1\le i\le N_{Q}\\1\le j\le N_{KV}}},\]
where, for any set $A$, $\mathbf 1_A$ is the associated indicator function associated to $A$.
In causal self-attention, the principle is that as the sentences are formed token by token in a sequential manner and the previous tokens are fixed, the addition of new tokens should not affect the previous parts of the sentence.
Hence, the query should only see itself and the tokens that have been defined before itself, i.e. $S_i=\{j\in\NN\,,\,j\le i\}$ and thus 
\[A=\left(\mathbf 1_{j\le i}\kappa\left({Q[i,:],K[j,:]}\right)\right)_{\substack{1\le i\le N_{Q}\\1\le j\le N_{KV}}}.\]
\begin{tcolorbox}
Equivalently, given the contextual mask $M$ such that 
\[ M_{ij}=\begin{cases}
    0&\text{if $j\in S_{i}$,}\\
    \varnothing_{\kappa}&\text{otherwise,}
\end{cases}\] 
(where $\varnothing_{\kappa}=-\infty$ for the exponential kernel). In the case of softmax attention, masked terms are given the value $-\infty$ so as to be attributed a zero probability after the exponential kernel function. 
Masked attention can then be seen as adding $M$ to the $QK^\T$ matrix before applying the kernel function giving 
$$A=\kappa( QK^{\T}+M).$$
Equivalently, the mask can be applied after the kernel by setting the masked elements to zero. By setting $\widetilde M_{ij}=\mathbf 1_{S_i}(j),$, with $\mathbf 1_{S_i}$, the indicating function for the set $S_i$, the masked attention is obtained as
$$A=\kappa( QK^{\T})\odot \widetilde M.$$
 \end{tcolorbox}
\paragraph{Cross-attention}
The second architectural difference in the decoder layers is that they use the encoder's output, $E_{N_L}$ in Figure \ref{fig:t_dec}, as the input of a second attention sublayer, which allows the decoder to relate its output sentence to the information contained in the original input. To do so, the intermediate state is used as the key and value inputs, $X^K$ and $X^V$, of a multi-headed attention sublayer as shown in  Figure \ref{fig:t_dec}. 
\begin{tcolorbox}
    \textbf{Remark:}\\
    Cross-attention \cite{alayracFlamingoVisualLanguage2022,rombachHighResolutionImageSynthesis2022a}, or more generally, shared embedding spaces \cite{kimViLTVisionandLanguageTransformer2021,radfordLearningTransferableVisual2021} are at the root of modern multimodal LLMs in which embedding subspaces can be shared through attention, or directly through a scalar product in the embedding space, between texts and images. 
\end{tcolorbox}
\subsection{The Transformer and its variants}
\begin{figure}
    \centering
    \includegraphics[width=\linewidth]{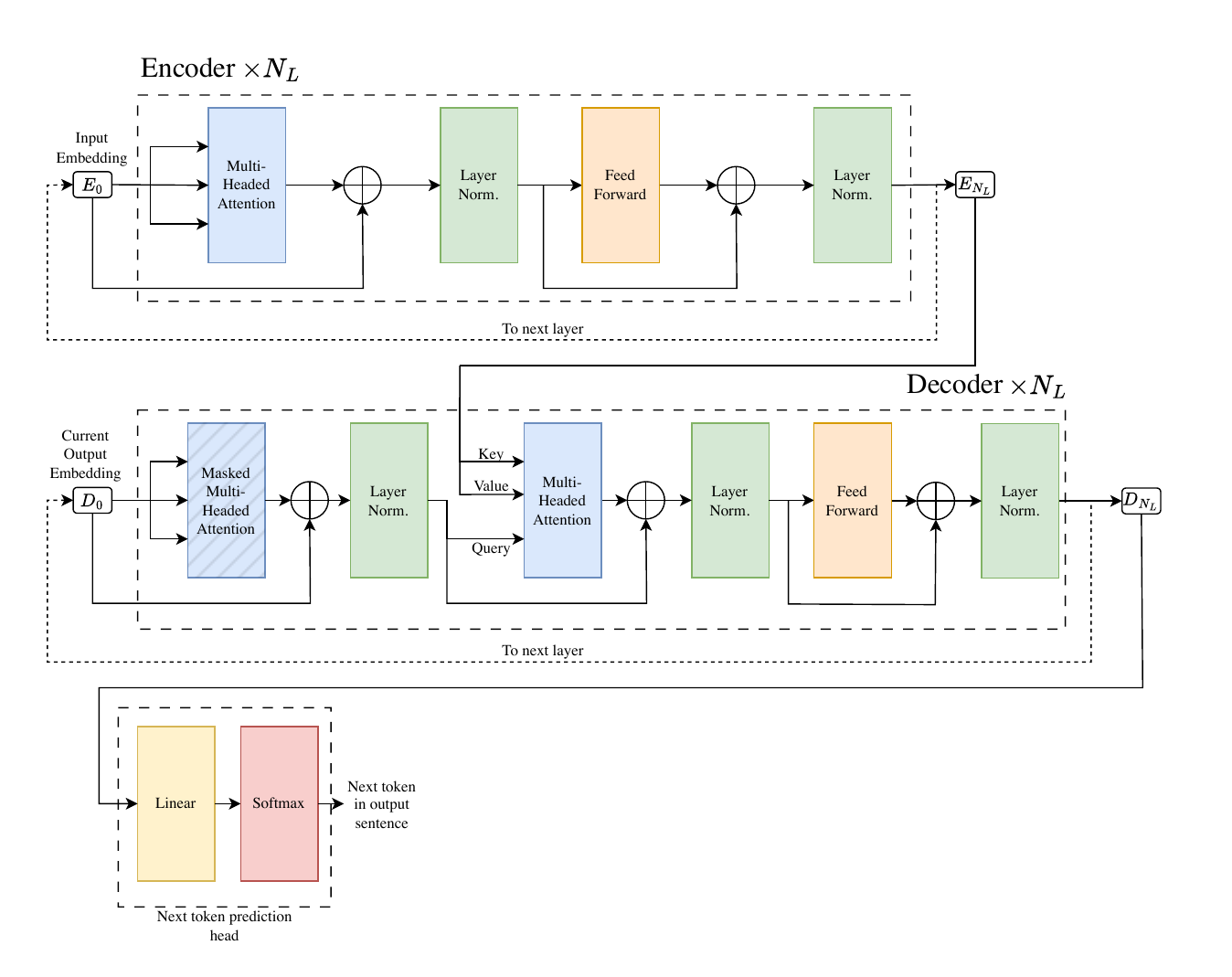}
    \caption{Transformer architecture.}
    \label{fig:t_arch}
\end{figure}
\begin{figure}
    \centering
    \includegraphics[width=\linewidth]{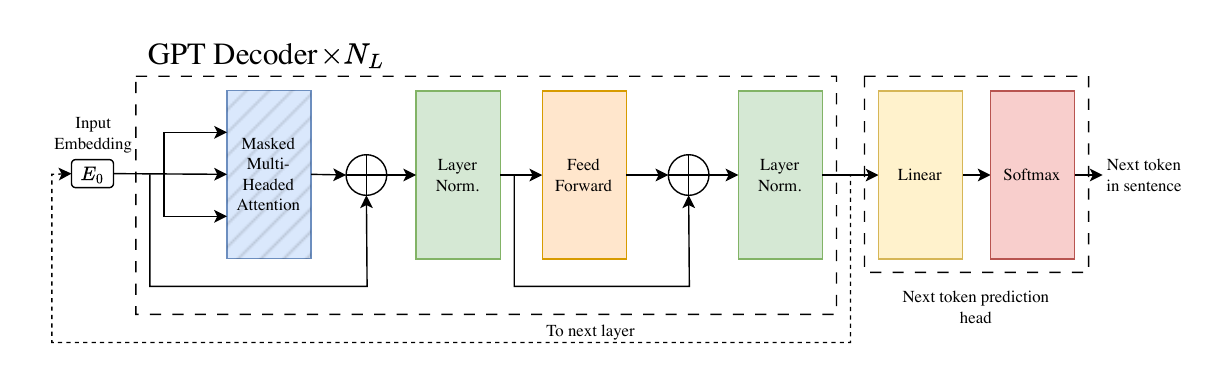}
    \caption{GPT architecture.}
    \label{fig:gpt_arch}
\end{figure}

The complete Transformer architecture, as can be seen in Figure~\ref{fig:t_arch}, is modelled for, and trained on, the task of machine translation and is at the origin of most LLMs today, albeit with several modifications. Indeed, while the Transformer network is well suited for sequence to sequence translation, a change in architecture is required for less structured problems.
Traditionally, these changes, which vary according to the task at hand, can be classified into two main groups, Encoder-only architectures such as BERT (Bidirectional Encoder Representation from Transformers) \cite{devlinBERTPretrainingDeep2019,heDeBERTaDecodingenhancedBERT2021,warnerSmarterBetterFaster2024}, more suited for tasks requiring the extraction of information from text such as token or sentence classification and Decoder-only architectures such as GPT(Generative Pre-Trained Transformers) \cite{radfordImprovingLanguageUnderstanding,brownLanguageModelsAre2020,touvronLLaMAOpenEfficient2023}, suited to generative tasks such as next-token prediction for text generation. The current trend, put forth in \cite{radfordImprovingLanguageUnderstanding}, lies in pretraining such models on tasks that require an understanding of the inherent structure of the text and with a huge amounts of available data, such as next-token prediction for the GPT architecture, using the architecture shown in Figure~\ref{fig:gpt_arch}, and masked text completion for Bert, and using the pretrained weights as an initialization for a fine-tuning a model with the same weights but a different output head for other tasks using smaller datasets.

\section{KV Caching, compression and attention optimization}
One of the main memory bottlenecks of modern LLMs arises due to what is known as \textbf{KV caching}. In the attention mechanism, the computation of the next output for the last token depends on the key and value vectors of all tokens. Indeed, this property is evident if we look at the last output of single-headed attention layer$$Y[N_{Q},j]=\left(Z^{-1}AV\right)[N_{Q},j]=\sum\limits_{i=1}^{N_{KV}}\frac{\kappa\left(Q[N_{Q},:],K[i,:]\right)}{\sum\limits_{\ell=1}^{N_{KV}}\kappa\left(Q[N_{Q},:],K[\ell,:]\right)}V[i,j],$$
 where the $N_Q$ index corresponds to the last query token. 
Note that the computational cost of the attention mechanism is of the order of $O(N_{Q}N_{KV}(d_{QK}+d_V))$ and that $d_{QK}$ and $d_V$ are generally small with respect to $N_Q$ and $N_{KV}$.
In the context of a chat for example, the token embeddings for the query, key and value terms are calculated in a stream, first you have the first sentence sent, then the answer, second sentence, etc... 
As the conversation goes on, the query, key and value vectors can either be recalculated at each new sentence, but with a computational cost of $O(N_{\text{tokens}}^3d)$, where $N_{KV}=N_{Q}=N_{\text{tokens}}$ as all the vectors need to be recomputed each time a new token appears.
Alternatively, the key and value vectors can be stored after each call to be used again for the next step, as we have seen above, once the key and value vectors are calculated, the query terms can be discarded as their only impact on subsequent predictions is through the KV vectors. 
KV caching consists in keeping these key and value vectors for each token in memory. This allows to obtain a linear cost in the number of tokens for each new token added( thus overall an $O(N_{\text{tokens}}^2d)$ complexity in total but an $O(N_QN_{\text{tokens}}d)$ cost for each new sentence of size $N_Q$). This does not come for free as the vectors need to be stored. As each we have to cache the KV vectors of each attention head for each layer the memory cost of such caching is $2N_{L}N_{h}N_{KV}d$ times the number of bits used per float. This can quickly become prohibitive, especially in the long-context case.
\begin{tcolorbox}
We call \textbf{streaming attention} the calculation of attention where the query input tokens are a subset of the key-value input tokens, as opposed to \textbf{full attention}. Note that the key-value vectors need to be cached for this to work as the previous query tokens are required for their construction.
\end{tcolorbox}
\subsection{Grouped Query Attention}
The size of the cache in memory is a bottleneck both in terms of storage and of communication overheads, methods have thus been developed to palliate to this.
One of the simplest, and most ubiquitous, such method consists in sharing the same key-value heads for multiple query heads, hence reducing the size of the vector that needs to be kept in memory. This is known as Grouped Query Attention (GQA), and in the limit of a single key-value head for all query heads it is sometimes referred to as Multi-Query Attention (MQA).  

One recent method which generalizes the use of lower-dimensional embedding spaces is the Multi-Headed Latent Attention mechanism developed by DeepSeek \cite{deepseek-aiDeepSeekV2StrongEconomical2024}. Let us now look at how this fits into our framework.

\subsection{Latent Attention}
The main idea behind Latent Attention is to construct a shared low-rank latent embedding from which the key and value vectors are formed for each head through a linear application. In the following, we set $X^{K}=X^{V}=X\in \RR^{N_{KV}\times d_{\text{in}}}$ and $d_{QK}=d_V=\dhead$.
Before diving into the formalism of Multi-headed Latent Attention(MLA), let us look back at Multi-Headed Attention and notice that the per-head embedding calculations can be seen as a single matrix multiplication,
\begin{equation*}
    \begin{aligned}
    \mathbf Q
    &=\begin{pmatrix}Q_{1}&\cdots& Q_{\Nh}\end{pmatrix}
    &=X_Q\begin{pmatrix}W^{Q}_{1}&\cdots & W^Q_{\Nh}\end{pmatrix}&\coloneq X_QW^{\mathbf Q} \in \RR^{N_{Q} \times \Nh \dhead} \\
    \mathbf K
    &=\begin{pmatrix}K_{1}&\cdots& K_{\Nh}\end{pmatrix}
    &=X\begin{pmatrix}W^{K}_{1}&\cdots & W^K_{\Nh}\end{pmatrix}
    &\coloneq XW^{\mathbf K} \in \RR^{N_{KV} \times \Nh \dhead}\\
    \mathbf V
    &=\begin{pmatrix}V_{1}&\cdots& V_{\Nh}\end{pmatrix}
    &=X\begin{pmatrix}W^{V}_{1}&\cdots & W^V_{\Nh}\end{pmatrix}
    &\coloneq XW^{\mathbf V}\in \RR^{N_{KV} \times \Nh \dhead}
    \end{aligned},
\end{equation*}
where $\mathbf Q$, $\mathbf K$, and $\mathbf V$ denote the column-wise concatenations of the head-specific query, key, and value matrices.
Similarly, we set
\begin{equation*}
W^{\mathbf Q}= \begin{pmatrix}
        W^{Q}_{1}  \cdots & W^{Q}_{\Nh}
       \end{pmatrix}\in\RR^{d_{in}\times\Nh\dhead},\quad 
W^{\mathbf K}
    = \begin{pmatrix}
        W^{K}_{1} & \cdots & W^{K}_{\Nh}
    \end{pmatrix}\in\RR^{d_{in}\times\Nh\dhead}\end{equation*}
       \begin{equation*}
 \text{ and }
W^{\mathbf V}
    = \begin{pmatrix}
        W^{V}_{1} & \cdots & W^{V}_{\Nh}
       \end{pmatrix}\in\RR^{d_{in}\times\Nh\dhead}.
\end{equation*} to denote the column-wise concatenations of the corresponding weight matrices.

Using this formulation, the latent attention method uses two latent subspaces, of dimensions $d_{L_Q}$ and $d_{L}\le d_{in}$, to construct a low-rank factorization of the query embeddings and a shared low-rank factorization of the key and value embeddings: 
\begin{equation*}
    \begin{aligned}
        W^{\mathbf Q}&=W^{L_Q}W^{\mathbf {L_Q Q}}=W^{L_Q}\begin{pmatrix}W^{{L_Q}Q}_{1}&\cdots & W^{{L_QQ}}_{\Nh}\end{pmatrix},\\
        W^{\mathbf K}&=W^{L}W^{\mathbf {LK}}=W^L\begin{pmatrix}W^{{LK}}_{1}&\cdots & W^{{LK}}_{\Nh}\end{pmatrix},\\
        W^{\mathbf V}&=W^{L}W^{\mathbf {LV}}=W^L\begin{pmatrix}W^{{LV}}_{1}&\cdots & W^{{LV}}_{\Nh}\end{pmatrix},
    \end{aligned}\label{eq:latent_KV}
\end{equation*}
where $W^{L}\in \RR^{d_{\text{in}}\times d_{L}}$ is the shared latent subspace weight matrix for the key and value embeddings, and we define for any $1\le h\le \Nh$, the matrices $W^{LK}_h\in \RR^{d_{L}\times \dhead}$, and $W_h^{LV}\in \RR^{d_{L}\times \dhead}$ to be the matrices that send the latent vectors onto the key and value embeddings of the $h$-th head. The matrices $W^{\mathbf{LK}}\in \RR^{d_{L}\times  \Nh \dhead}$ and $W^{\mathbf{LV}}\in \RR^{d_{L}\times  \Nh \dhead}$ are the the column-wise concatenation of $\Big(W^{LK}_h\Big)_{1\le h\le \Nh}$, and respectively $\Big(W^{LV}_h\Big)_{1\le h\le \Nh}$. Analogously, $W^{L_Q}\in \RR^{d_{\text{in}}\times d_{L_Q}}$ is latent subspace weight matrix for the query embeddings and $W^{\mathbf{L_QQ}}$ is the column-wise concatenation of $\Big(W^{L_QQ}_h\Big)_{1\le h\le \Nh}$, the matrices that send the latent query embeddings back to the shared query-key subspace.
\begin{tcolorbox}
    In practice, in the DeepSeek V2 model \cite{deepseek-aiDeepSeekV2StrongEconomical2024}, the latent model is trained directly and is not created from the key and value matrices using the above factorizations. 
\end{tcolorbox}

This idea allows to keep a single cache vector per token, \textit{shared between all heads}; the latent embedding $L=XW^L$, instead of $K$ and $V$. Furthermore, the latent formulation allows to merge weight matrices by rewriting the attention mechanism.
First, for each head $1\le h\le \Nh$, the Latent-to-Query weights $W_h^{L_QQ}$ and the Latent-to-Key weights $W_h^{LK}$ can be merged, reducing the computational cost as $d_L$ and $d_{L_Q}<d_{in}$. Indeed, we can rewrite, using \eqref{eq:latent_KV}, 
$$ Q_{h}K_{h}^\T=L_QW^{L_QQ}_{h}{W^{LK}_{h}}^{\T}{W^{L}}^{\T}X^{\T}=L_Q\underbrace{W^{L_QQ}_{h}{W^{LK}_{h}}^{\T}}_{{W^{LQK}_{h}}\in\RR^{d_{L_Q}\times d_{L}}}L^{\T},$$
which gives us, by indexing on each head $h$, the merged query tensor $\left({W^{LQK}_{h}}\right)_{1\le h\le\Nh}\in\RR^{d_{L_Q}\times d_{L}\times\Nh}$.

Second, the Latent-to-Value weights $W^{LV}$ can be merged with the head combining output weights $W^O$.
Indeed, we can rewrite the output of the attention mechanism as follows 
\[
Y=\operatorname{concat}_{\mathrm{col}}\left(Z_{h}^{-1}A_{h}V_{h}\right)_{1\le h\le \Nh}W^{O}=\underbrace{\operatorname{concat}_{\mathrm{col}}\left(Z_{h}^{-1}A_{h}L\right)_{1\le h\le \Nh}}_{Y '}\underbrace{\begin{pmatrix}W^{LV}_{1}&&0\\&\ddots&\\0&&W^{LV}_{\Nh}\end{pmatrix}W^{O}}_{{W^{LO}}},
\]
where we highlight again the fact that $L$ does not depend on $h$ as it is shared for all heads.
This gives us the new latent output matrix $W^{LO}\in\RR^{\Nh d_L \times d_{out}}$.

The latent attention mechanism can thus be simplified to 
$$Y= \operatorname{concat}_{\mathrm{col}}\left(Z_{h}^{-1}A_{h}L\right)_{1\le h\le \Nh}{W^{LO}},$$with $A_{h}=\exp(L_QW_h^{LQK}L^{\T})$ and $Z_h=\diag\left(\sum\limits_{j=1}^{N_{KV}} \exp\Big( (L_QW^{LQK}_hL^\T)[i,j]\Big)\right)_{1\le i\le N_{Q}}$. 
With this formulation, note that the only model weights that need to be kept are $W^{LO}$, $W^L$, $W^{L_Q}$ and ${W^{LQK}}$ and only the latent vectors need to be kept in memory.
\begin{tcolorbox}
    In the absence of positional embeddings, latent attention models can be represented exactly in the form of an equivalent GQA/MHA model thanks to the low-rank factorizations shown above. However, in practice, when positional embeddings are applied, this equivalence is lost.
    This loss of equivalence is due to the fact that positional embeddings, such as the Rotary Position Embeddings \cite{suRoFormerEnhancedTransformer2023} (RoPE), in classical attention are applied after constructing $K$. For RoPE, in latent attention, this would amount to applying a position-dependent matrix $R_m\in\RR^{\dhead\times \dhead}$, as follows $$Q_h[i,:]R_iR_j^\T K_h[j,:]^\T=(L_Q[i,:]) W_h^{L_QQ} R_i R_j^\T {W_h^{LK}}^\T (L[j,:])^\T.$$
    This would not allow the merging of the $W^{LK}_h$ and $W^{L_Q Q}_h$ matrices and would require recomputing the application of the positional embedding at each evaluation, not only preventing a speedup but incurring an additional computational overhead. The latent attention RoPE variant  instead consist in appending to the key and query vectors ``non-latent'' part on which the positional embedding is applied, breaking the equivalence but keeping the computational advantage of the latent model.
    Recently, methods to approximately convert pretrained GQA/MHA models of attention, from previously learnt models, to the Latent attention formalism have been developed \cite{mengTransMLAMultiHeadLatent2025} and allow for the use of the DeepSeek's optimizations on other models. 
    
\end{tcolorbox}
The tables \ref{tab:mha_mem} and \ref{tab:mla_mem}  summarize the matrices and caches that need to be stored in memory for the MHA and MLA cases respectively. Example of model sizes are given in table \ref{tab:llama3_vs_gemma3_vs_deepseekV2}.

\begin{table}[h!]
    \centering
    \begin{minipage}[t]{0.48\textwidth} 
        \centering
        \begin{tabular}{|l|l|l|}
            \hline
            \textbf{Tensor/Matrix} & \textbf{Dimensions} & \textbf{Description} \\ \hline
            $K$ & $\RR^{N_{KV}\times d_{QK}\times\Nh}$ & \multirow{2}{*}{KV cache} \\ \cline{1-2}
            $V$ & $\RR^{N_{KV}\times \dhead\times\Nh}$ & \\ \hline
            $W^{Q}$ & $\RR^{d_{\text{in}}\times d_{QK} \times \Nh}$ & \multirow{4}{*}{Weights} \\ \cline{1-2}
            $W^{K}$ & $\RR^{d_{\text{in}}  \times d_{QK}\times \Nh}$ & \\ \cline{1-2}
            $W^{V}$ & $\RR^{d_{\text{in}} \times \dhead\times \Nh}$ & \\ \cline{1-2}
            $W^O$ & $\RR^{\Nh \dhead \times d_{\text{out}}}$ & \\ \hline
        \end{tabular}
        \vspace{0.5em}
        \captionof{table}{Dimensions of tensors and matrices stored in multi-headed attention. The total number of floats to be stored is $ \Nh\ N_{KV} \cdot \left(d_{QK}+\dhead\right)  +2\Nh d_{\text{in}} d_{QK} +\Nh\dhead \left(d_{\text{in}}+d_{\text{out}}\right)$.
        For GQA models, the $\Nh$ in the dimensions needs to be replaced with the number of KV heads for $W^K$,$W^V$ and the cache terms.}
        \label{tab:mha_mem}
    \end{minipage}
    \hfill 
    \begin{minipage}[t]{0.48\textwidth} 
        \centering
        \begin{tabular}{|l|l|l|}
            \hline
            \textbf{Tensor/Matrix} & \textbf{Dimensions} & \textbf{Description} \\ \hline
            $L$ & $\RR^{N_{KV} \times d_L}$ & L cache \\ \hline
            $W^L$ & $\RR^{d_{\text{in}}\times d_{L}}$ & \multirow{3}{*}{Weights} \\ \cline{1-2}
            ${W^{ Q}}'$ & $\RR^{d_{\text{in}}\times d_{L}\times \Nh}$ & \\ \cline{1-2}
            ${W^O}'$ & $\RR^{\Nh d_L \times d_{\text{out}}}$ & \\ \hline
        \end{tabular}
                \vspace{0.5em}
        \captionof{table}{Dimensions of tensors and matrices stored in multi-headed  latent attention. The total number of floats to be stored is $(N_{KV} \cdot d_{L}) + (d_{\text{in}} \cdot d_L) + (\Nh \cdot d_{L} \cdot (d_{\text{in}}+d_{\text{out}}))$.}
        \label{tab:mla_mem}
    \end{minipage}
\end{table}

\begin{table}[h!]
\centering
\begin{tabular}{|l|c|c|c|}
\hline
\textbf{Parameter} & \textbf{Llama 3 70B} & \textbf{Gemma 3 27B}& \textbf{Deepseek V2}\\
\hline
Number of Layers & 80 & 62 &  60\\
\hline
Number of Heads $\Nh$& 64 & 32 & 128 \\
\hline
Number of Key/Value Heads & 8 & 16 & - \\
\hline
Hidden Dimensions $d_{\text{in}}=d_\text{out}$ & 8,192 & 5,376 & 5120 \\
\hline
KV Dimension $\dhead=d_{QK}$ & 128 & 128 & 512($d_L$) \\
\hline
\end{tabular}
\vspace{0.5em}
\caption{Comparison of Llama 3 70B \cite{grattafioriLlama3Herd2024} and Gemma 3 27B \cite{teamGemma3Technical2025} and DeepSeek V2 Model \cite{deepseek-aiDeepSeekV2StrongEconomical2024} Specifications. Note that DeepSeekV2 is a mixture-of-experts model with latent attention while Gemma 3 and Llama 3 are Grouped Query attention models. These models use a Grouped Query attention model wherein multiple attention mechanisms are used in parallel, these are known as attention heads. To reduce the computational overhead, the key and value matrices are shared between heads, the number of key/value matrices is denoted as the number of KV heads above.
}
\label{tab:llama3_vs_gemma3_vs_deepseekV2}
\end{table}
\pagebreak
\section{Acknowledgments}
The author would like to thank Alice Cortinovis and Laura Grigori for their valuable feedback on this document.
\printbibliography
\end{document}